\documentclass[a4paper]{amsart}

\RequirePackage{amsmath}
\RequirePackage{bm}
\RequirePackage{amssymb}
\RequirePackage{upref}
\RequirePackage{amsthm}
\RequirePackage{enumerate}
\RequirePackage{pb-diagram}
\RequirePackage{amsfonts}
\RequirePackage[mathscr]{eucal}
\RequirePackage{verbatim}
\RequirePackage{xr}
\RequirePackage{graphicx}
\usepackage{calc}
\usepackage{xspace}
\RequirePackage{color}
\RequirePackage{ifthen}

\newcommand{\cf}{cf.\@\xspace}
\newcommand{\resp}{resp.\@\xspace}

\newcommand{\al}{\alpha}
\newcommand{\bet}{\beta}
\newcommand{\ga}{\gamma}
\newcommand{\de}{\delta }
\newcommand{\e}{\epsilon}

\newcommand{\f}{\varphi}
\newcommand{\h}{\eta}

\newcommand{\ka}{\kappa}
\newcommand{\lam}{\lambda}

\newcommand{\C}{\varGamma}
\newcommand{\D}{\varDelta}

\newcommand{\Lam}{\varLambda}

\newcommand{\msp[1]}[1]{\mspace{#1mu}}

\newcommand{\R}[1][n+1]{{\protect\mathbb R}^{#1}}

\newcommand{\N}{{\protect\mathbb N}}

\newcommand{\eR}{\stackrel{\lower1ex \hbox{\rule{6.5pt}{0.5pt}}}{\msp[3]\R[]}}
\newcommand{\eN}{\stackrel{\lower1ex \hbox{\rule{6.5pt}{0.5pt}}}{\msp[1]\N}}
\newcommand{\eO}{\stackrel{\lower1ex \hbox{\rule{6pt}{0.5pt}}}{\msc O}}
\newcommand{\mf}[1]{\mathfrak {#1}}

\DeclareMathOperator{\ad}{ad}
\DeclareMathOperator{\Ad}{Ad}

\DeclareMathOperator{\tr}{tr}
\DeclareMathOperator{\SO}{SO}

\newcommand\ra{\rightarrow}

\newcommand{\un}{\infty}
\newcommand{\A}{\forall}

\newcommand{\uu}{\cup}

\newcommand{\uuu}{\bigcup}

\newcommand{\uud}{ \stackrel{\lower 1ex \hbox {.}}{\uu}}
\newcommand{\uuud}[1]{ \stackrel{\lower 1ex \hbox {.}}{\uuu_{#1}}}

\newcommand{\sminus}[1][28]{\raise 0.#1ex\hbox{$\scriptstyle\setminus$}}

\newcommand{\abs}[1]{\lvert#1\rvert}

\newcommand{\riema}[4]{{\bar R}_{#1#2#3#4}}

\newcommand{\tit}{\textit}

\newcommand{\tup}{\textup}

\newcommand{\mc}{\protect\mathcal}
\newcommand{\msc}{\protect\mathscr}

\providecommand{\bysame}{\makebox[3em]{\hrulefill}\thinspace}

\newcommand{\bt}{\begin{thm}}
\newcommand{\bl}{\begin{lem}}
\newcommand{\bc}{\begin{cor}}
\newcommand{\bd}{\begin{definition}}
\newcommand{\bpp}{\begin{prop}}
\newcommand{\br}{\begin{rem}}
\newcommand{\bn}{\begin{note}}
\newcommand{\be}{\begin{ex}}
\newcommand{\bes}{\begin{exs}}
\newcommand{\bb}{\begin{example}}
\newcommand{\bbs}{\begin{examples}}
\newcommand{\ba}{\begin{axiom}}
\newcommand{\bas}{\begin{assumption}}

\newcommand{\et}{\end{thm}}
\newcommand{\el}{\end{lem}}
\newcommand{\ec}{\end{cor}}
\newcommand{\ed}{\end{definition}}
\newcommand{\epp}{\end{prop}}
\newcommand{\er}{\end{rem}}
\newcommand{\en}{\end{note}}
\newcommand{\ee}{\end{ex}}
\newcommand{\ees}{\end{exs}}
\newcommand{\eb}{\end{example}}
\newcommand{\ebs}{\end{examples}}
\newcommand{\ea}{\end{axiom}}
\newcommand{\eas}{\end{assumption}}

\newcommand{\bp}{\begin{proof}}
\newcommand{\ep}{\end{proof}}
\newcommand{\eps}{\renewcommand{\qed}{}\end{proof}}

\newcommand{\bal}{\begin{align}}

\newcommand{\bi}[1][1.]{\begin{enumerate}[\upshape #1]}
\newcommand{\bia}[1][(1)]{\begin{enumerate}[\upshape #1]}
\newcommand{\bin}[1][1]{\begin{enumerate}[\upshape\bfseries #1]}
\newcommand{\bir}[1][(i)]{\begin{enumerate}[\upshape #1]}
\newcommand{\bic}[1][(i)]{\begin{enumerate}[\upshape\hspace{2\cma}#1]}
\newcommand{\bis}[2][1.]{\begin{enumerate}[\upshape\hspace{#2\parindent}#1]}
\newcommand{\ei}{\end{enumerate}}

\newcommand\ndots{\raise 0.47ex \hbox {,}\hskip0.06em\cdots %
     \raise 0.47ex \hbox {,}\hskip0.06em} 

\newcommand{\q}{\quad}
\newcommand{\qq}{\qquad}

\newcommand{\hp}{\hphantom}

\newcommand\nd{\noindent}

\newskip\Csmallskipamount                                                
\Csmallskipamount=\smallskipamount
\newskip\Cmedskipamount
\Cmedskipamount=\medskipamount
\newskip\Cbigskipamount
\Cbigskipamount=\bigskipamount

\newcommand\cvs{\vspace\Csmallskipamount}   
\newcommand\cvm{\vspace\Cmedskipamount}

\newskip\csa
\csa=\smallskipamount

\newskip\cma
\cma=\medskipamount

\newskip\cba
\cba=\bigskipamount

\newdimen\spt
\newcommand\citem{\cvs\advance\itemno by
1{(\romannumeral\the\itemno})\hskip3pt}
\newcommand{\bitem}{\cvm\nd\advance\itemno by
1{\bf\the\itemno}\hspace{\cma}}

\newcount\itemno
\newcommand{\lae}[1]{\label{E:#1}}
\newcommand{\lat}[1]{\label{T:#1}}

\newcommand{\re}[1]{\eqref{E:#1}}
\newcommand{\frt}[1]{Theorem~\ref{T:#1} on page~\tup{\pageref{T:#1}}}

\newcommand{\fre}[1]{\eqref{E:#1} on page~\tup{\pageref{E:#1}}}

\newskip\thmskip
\thmskip=\parindent

\newskip\hsk
\newenvironment{hinw}{\labelsep=0pt\begin{list}{}{\labelsep=0pt\itemindent=0pt\labelwidth=0pt\leftmargin=\parindent\rightmargin=0pt\partopsep=\cba}%
\item\it\nopagebreak\nopagebreak}%
{\end{list}}

\newcommand\bh{\begin{hinw}}
\newcommand{\eh}{\end{hinw}}

\newtheoremstyle{normal}
  {\cba}
  {\cba}
  {}
  {\thmskip}
  {\bfseries}
  {.}
  {\hsk}
  {}

\newtheoremstyle{abschnitt}
  {\cba}
  {\cba}
  {}
  {\thmskip}
  {\bfseries}
  {.}
  {\hsk}
  {}

\newtheoremstyle{italic}
  {\cba}
  {\cba}
  {\itshape}
  {\thmskip}
  {\bfseries}
  {.}
  {\hsk}
  {}

\newtheoremstyle{aufgaben}
  {\cba}
  {\cba}
  {}
  {}
  {\normalsize\bfseries}
  {.}
  {\hsk}
  {}

\newtheoremstyle{break}
  {\cba}
  {\cba}
  {\itshape}
  {}
  {\bfseries}
  {.}
  {\newline}
  {}

\swapnumbers
\theoremstyle{italic}
\newtheorem{thm}[subsection]{Theorem}
\newtheorem{lem}[subsection]{Lemma}
\newtheorem{prop}[subsection]{Proposition}
\newtheorem{cor}[subsection]{Corollary}

\theoremstyle{normal}
\newtheorem{rem}[subsection]{Remark}
\newtheorem{definition}[subsection]{Definition}
\newtheorem{example}[subsection]{Example}
\newtheorem{examples}[subsection]{Examples}
\newtheorem{ex}[subsection]{Exercise}
\newtheorem{note}[subsection]{}
\newtheorem{axiom}[subsection]{Axiom}
\newtheorem{assumption}[subsection]{Assumption}

\theoremstyle{aufgaben}
\newtheorem{exs}[subsection]{Exercises}

\swapnumbers

\numberwithin{equation}{section}
\numberwithin{figure}{section}

\newenvironment{textequation}[1][0.8]
{\begin{equation}
\begin{aligned}
\begin{minipage}{#1\linewidth}}
{\end{minipage}
\end{aligned}
\end{equation}
\ignorespacesafterend}

\newcommand{\btext}{\begin{textequation}}
\newcommand{\etext}{\end{textequation}}

\def\hinweis{\@startsection{subsection}{2}%
 \z@{0.7\linespacing\@plus 0.5\linespacing}{0.7\linespacing}%
{\normalfont\itshape\indent}}

\newcounter{hours}\newcounter{minutes}
\newcommand{\printtime}{%
\setcounter{hours}{\time/60}%
\setcounter{minutes}{\time-\value{hours}*60}%
\ifthenelse{\value{minutes}<9}{\thehours :0\theminutes}{\thehours:\theminutes}}

\usepackage[german,english]{babel}
\usepackage{graphicx}
\RequirePackage{amsmath}
\RequirePackage{bm}
\RequirePackage{amssymb}
\RequirePackage{upref}
\RequirePackage{amsthm}
\RequirePackage{enumerate}
\RequirePackage{pb-diagram}
\RequirePackage{amsfonts}
\RequirePackage[mathscr]{eucal}
\RequirePackage{verbatim}
\RequirePackage{xr}
\RequirePackage{graphicx}
\usepackage{calc}
\usepackage{xspace}

\makeatletter
\RequirePackage{color}
\newcommand{\ann}[1]{\renewcommand{\@makefnmark}{\mbox{$^{\color{red}{\@thefnmark}}$}}%
\footnote {#1}}
\makeatother

\RequirePackage{upref}
\RequirePackage{amsthm}
\RequirePackage{enumerate}
\usepackage[mathscr]{eucal}

\usepackage{xr-hyper}

\usepackage{calc}

\newlength{\oddsidemarginlength}
\newlength{\topmarginlength}

\newcounter{numberoflines}
\newcounter{tempcc}
\oddsidemargin=\oddsidemarginlength
\evensidemargin=\oddsidemargin
\topmargin=\topmarginlength
\begin{document}

\flushbottom


\title{An energy gap for Yang-Mills connections}

\author{Claus Gerhardt}
\address{Ruprecht-Karls-Universit\"at, Institut f\"ur Angewandte Mathematik,
Im Neuenheimer Feld 294, 69120 Heidelberg, Germany}
\email{gerhardt@math.uni-heidelberg.de}
\urladdr{http://www.math.uni-heidelberg.de/studinfo/gerhardt/}
\thanks{This work has been supported by the DFG}

%
\date{\today}
%


\begin{abstract}
Consider a Yang-Mills connection over a Riemann manifold $M=M^n$, $n\ge 3$, where $M$ may be compact or complete. Then its energy must be bounded from below by some positive constant, if $M$ satisfies certain conditions, unless the connection is flat.
\end{abstract}

\maketitle

\tableofcontents

\setcounter{section}{0}
\section{Introduction}
We consider the problem: When is a Yang-Mills connection non-flat? Of course, the trivial answer $F_{\mu\lam}\not\equiv 0$ is unsatisfactory. Bourguignon and Lawson proved in \cite[Theorem C]{bourguignon}, among other results, that any Yang-Mills connection over $S^n$, $n\ge 3$, the field strength of which satisfies the pointwise estimate
\begin{equation}
F^2=-\tr(F_{\mu\lam}F^{\mu\lam})<\binom n2
\end{equation}
is flat.

We want to prove that under certain assumptions on the base space $M$, which is supposed to be a Riemannian manifold of dimension $n\ge 3$, the energy of a Yang-Mills connection has to satisfy
\begin{equation}\lae{1.2}
\Big(\int_M\abs F^\frac n2\Big)^\frac2n\ge \ka_0>0,
\end{equation}
where $\ka_0$ depends only on the Sobolev constants of $M$, $n$ and the dimension of the Lie group $\mc G$, unless the connection  is flat.

Here,
\begin{equation}
\abs F= \sqrt{F^2},
\end{equation}
and we also call the left-hand side of \re{1.2} energy though this label is only correct when $n=4$. However, this norm is also  the crucial norm, which has to be (locally) small, used to prove  regularity of a connection, \cf \cite[Theorem 1.3]{uhlenbeck:lp}. 

The exponent $\frac n2$ naturally pops up when Sobolev inequalities are applied to solutions of differential equations satisfied by the field strength or the energy density of a connection in the adjoint bundle.

We distinguish two cases: $M$ compact and $M$ complete and non-compact. When $M$ is compact, we require 
\begin{equation}\lae{1.3}
\bar R_{\al\bet}\Lam^\al_\lam\Lam^{\bet\lam}-\tfrac12\riema\al\bet\mu\lam\Lam^{\al\bet}\Lam^{\mu\lam}\ge c_0\Lam_{\al\bet}\Lam^{\al\bet}
\end{equation}
for all skew-symmetric $\Lam_{\al\bet}\in T^{0,2}(M)$, where $0<c_0$, while for non-compact $M$ the weaker assumption
\begin{equation}\lae{1.4}
\bar R_{\al\bet}\Lam^\al_\lam\Lam^{\bet\lam}-\tfrac12\riema\al\bet\mu\lam\Lam^{\al\bet}\Lam^{\mu\lam}\ge0
\end{equation}
and in addition 
\begin{equation}\lae{1.6}
\Big(\int_Mu^{\frac{2n}{n-2}}\Big)^{\frac{n-2}n}\le c_1 \int_M\abs{Du}^2\qq\A\, u\in H^{1,2}(M)
\end{equation}
should be satisfied.

\br
(i)  If $M$ is a space of constant curvature
\begin{equation}\lae{1.7}
\riema\al\bet\mu\lam=K_M(\bar g_{\al\mu}\bar g_{\bet\lam}-\bar g_{\al\lam}\bar g_{\bet\mu}),
\end{equation}
then
\begin{equation}
\bar R_{\al\bet}\Lam^\al_\lam\Lam^{\bet\lam}-\tfrac12\riema\al\bet\mu\lam\Lam^{\al\bet}\Lam^{\mu\lam}=(n-2)K_M\Lam_{\al\bet}\Lam^{\al\bet}.
\end{equation}
In case $n=2$ the curvature term therefore vanishes, and this result is also valid for an arbitrary two-dimensional Riemannian manifold, since the curvature tensor then has the same structure as in \re{1.7} though $K_M$ is not necessarily constant.

(ii) If $M=\R[n]$, $n\ge 3$, the conditions \re{1.4} and \re{1.6} are always valid.
\er

\bt\lat{1.2} 
Let $M=M^n$, $n\ge 3$, be a compact Riemannian  for which the condition \re{1.3} with $c_0>0$ holds. Then any Yang-Mills connection over $M$ with compact, semi-simple Lie group is either flat or satisfies \re{1.2} for some constant $\ka_0>0$ depending on the Sobolev constants of $M$, $n,c_0$, and the dimension of the Lie group.
\et

\bt\lat{1.3}
Let $M=M^n$, $n\ge 3$, be complete, non-compact and assume that the conditions \re{1.4} and\re{1.6} hold. Then any Yang-Mills connection over $M$ with compact, semi-simple Lie group is either flat or the estimate \re{1.2} is valid. The constant $\ka_0>0$ in \re{1.2} depends on the constant $c_1$ in \re{1.6}, $n$, and the dimension of the Lie group.
\et

\section{The compact case}

Let $(P,M,\mc G,\mc G)$ be a principal fiber bundle where $M=M^n$, $n\ge 3$ is a compact Riemannian manifold with metric $\bar g_{\al\bet}$ and $\mc G$ a compact, semi-simple Lie group with Lie algebra $\mf g$. Let $f_c=(f^a_{cb})$ be a basis of $\ad\mf g$ and
\begin{equation}
A_\mu=f_cA^c_\mu
\end{equation}
a Yang-Mills connection in the adjoint bundle $(E,M,\mf g,\Ad(\mc G))$.

The curvature tensor of the connection is given by
\begin{equation}\lae{2.2}
R^a_{\hp{a} b\mu\lam}=f^a_{cb}F^c_{\mu\lam},
\end{equation}
where
\begin{equation}
F_{\mu\lam}=f_cF^c_{\mu\lam}
\end{equation}
is the field strength of the connection, and 
\begin{equation}
F^2\equiv \ga_{ab}F^a_{\mu\lam}F^{b\mu\lam}=R_{ab\mu\lam}R^{ab\mu\lam}
\end{equation}
the energy density of the connection---at least up to a factor $\frac14$.

Here, $\ga_{ab}$ is the Cartan-Killing metric acting on elements of the fiber $\mf g$, and Latin indices are raised or lowered with respect to the inverse $\ga^{ab}$ or $\ga_{ab}$, and Greek indices with respect to the metric of $M$.

\bd
The adjoint bundle $E$ is vector bundle; let $E^*$ be the dual bundle, then we denote by 
\begin{equation}
T^{r,s}(E)=\C(\underset{r}{\underbrace{E\otimes\cdots\otimes E}}\otimes\underset{s}{\underbrace{E^*\otimes\cdots\otimes E^*}})
\end{equation}
the sections of the corresponding tensor bundle.
\ed

Thus, we have
\begin{equation}
F^a_{\mu\lam}\in T^{1,0}(E)\otimes T^{0,2}(M).
\end{equation}

Since $A_\mu$ is a Yang-Mills connection it solves the Yang-Mills equation
\begin{equation}\lae{2.7}
F^{a\al}_{\hp{a\al}\lam;\al}=0,
\end{equation}
where we use Einstein's summation convention, a semi-colon indicates covariant differentiation, and where we stipulate that a covariant derivative is always a \tit{full} tensor, i.e.,
\begin{equation}
F^a_{\mu\lam;\al}=F^a_{\mu\lam,\al}+f^a_{bc}A^b_\al F^c_{\mu\lam}-\bar\C^\ga_{\al\mu} F^a_{\ga\lam}-\bar \C^\ga_{\al\lam}F^a_{\mu\ga},
\end{equation}
where $\bar\C^\ga_{\al\bet}$ are the Christoffel symbols of the Riemannian connection; a comma indicates partial differentiation.

Before we formulate the crucial lemma let us note that $\riema\al\bet\ga\de$ \resp $\bar R_{\al\bet}$ symbolize the Riemann curvature tensor \resp the Ricci tensor of $\bar g_{\al\bet}$.

\bl
Let $A_\mu$ be a Yang-Mills connection, then its energy density $F^2$ solves the equation
\begin{equation}\lae{2.9}
\begin{aligned}
&-\tfrac14\D F^2+\tfrac12 F_{a\mu\lam;\al}F^{a\mu\lam\hp{;}\al}_{\hp{a\mu\lam};}+\bar R_{\bet\mu}F^{a\bet}_{\hp{a\bet}\lam}F^{\hp a\mu\lam}_a-\tfrac12 \riema\al\bet\mu\lam F^{\hp a\al\bet}_aF^{a\mu\lam}\\
&=-f^a_{cb}F^c_{\al\mu}F^{b\al}_{\hp{b\al}\lam}F^{\hp a\mu\lam}_a.
\end{aligned}
\end{equation}
\el

\bp
Differentiating \re{2.7} covariantly with respect to $x^\mu$ and using the Ricci identities we obtain
\begin{equation}\lae{2.10}
\begin{aligned}
0&=-F^{a\al}_{\hp{a\al}\lam;\al\mu}\\
&=-F^{a\al}_{\hp{a\al}\lam;\mu\al}+R^a_{\hp a b\al\mu}F^{b\al}_{\hp{b\al}\lam}+\bar R^\al_{\hp\al\bet\al\mu} F^{a\bet}_{\hp{a\al}\lam} +\bar R^\bet_{\hp\al\lam\mu\al} F^{a\al}_{\hp{a\al}\bet}.
\end{aligned}
\end{equation}

On the other hand, differentiating the second Bianchi identities
\begin{equation}
0=F^a_{\al\lam;\mu}+F^a_{\mu\al;\lam}+F^a_{\lam\mu;\al}
\end{equation}
we infer
\begin{equation}
0=F^{a\al}_{\hp{a\al}\lam;\mu\al}+F^{a\hp{\mu}\al}_{\hp a\mu\hp\al ;\lam\al}+\D F^a_{\lam\mu},
\end{equation}
and we deduce further
\begin{equation}
-\D F^a_{\mu\lam}F_a^{\hp a\mu\lam}=-2F^{a\al}_{\hp{a\al}\lam;\mu\al}F_a^{\hp a\mu\lam}.
\end{equation}

In view of \re{2.10} we then conclude
\begin{equation}
\begin{aligned}
0=-\tfrac12\D F^a_{\mu\lam}F_a^{\hp a\mu\lam}&+R^a_{\hp a b\al\mu}F^{b\al}_{\hp{b\al}\lam}F_a^{\hp a\mu\lam} +\bar R_{\bet\mu}F^{a\bet}_{\hp{a\bet}\lam}F_a^{\hp a\mu\lam}\\
&+\bar R^\bet_{\hp\bet \lam\mu\al}F^{a\al}_{\hp{a\al}\bet}F_a^{\hp a\mu\lam},
\end{aligned}
\end{equation}
which is equivalent to
\begin{equation}
\begin{aligned}
0=-\tfrac12\D F^a_{\mu\lam}F_a^{\hp a\mu\lam}&+f^a_{cb}F^c_{\al\mu}F^{b\al}_{\hp{b\al}\lam}F_a^{\hp a\mu\lam} +\bar R_{\bet\mu}F^{a\bet}_{\hp{a\bet}\lam}F_a^{\hp a\mu\lam}\\
&-\bar R_{\al\mu\bet\lam}F^{a\al\bet}F_a^{\hp a\mu\lam},
\end{aligned}
\end{equation}
in view of \re{2.2}.

Finally, using the first Bianchi identities,
\begin{equation}
\riema\al\bet\mu\lam+\riema\al\mu\lam\bet+\riema\al\lam\bet\mu=0,
\end{equation}
we deduce
\begin{equation}
\begin{aligned}
\riema\al\bet\mu\lam  F^{a\al\bet}F_a^{\hp a\mu\lam}+\riema\al\mu\lam\bet F^{a\al\bet}F_a^{\hp a\mu\lam}+\riema\al\lam\bet\mu F^{a\al\bet}F_a^{\hp a\mu\lam}=0,
\end{aligned}
\end{equation}
and hence
\begin{equation}
\riema\al\bet\mu\lam F^{a\al\bet}F_a^{\hp a\mu\lam}=2 \riema \al\mu\bet\lam F^{a\al\bet}F_a^{\hp a\mu\lam},
\end{equation}
from which the equation \re{2.9} immediately follows.
\ep

\bp[Proof of \frt{1.2}] Define 
\begin{equation}
u=F^2,
\end{equation}
then
\begin{equation}
\bar R_{\bet\mu}F^{a\bet}_{\hp{a\bet}\lam}F^{\hp a\mu\lam}_a-\tfrac12 \riema\al\bet\mu\lam F^{\hp a\al\bet}_aF^{a\mu\lam}\ge c_0 u,
\end{equation}
where $c_0>0$, in view of the assumption \fre{1.3}.

Multiplying \re{2.9} with $u$ and integrating by part we obtain
\begin{equation}\lae{2.23}
\tfrac38\int_M\abs{Du}^2+c_0\int_Mu^2\le c\int_M\sqrt uu^2,
\end{equation}
where we used the simple estimate
\begin{equation}
\abs{Du}^2\le 4 F_{a\mu\lam;\al}F^{a\mu\lam\hp{;}\al}_{\hp{a\mu\lam};}u^2
\end{equation}
and where $c$ depends on $n$ and the dimension of $\mf g$; note that
\begin{equation}
f_c\in \SO(\mf g, \ga_{ab}).
\end{equation}

The integral on the right-hand side of \re{2.23} is estimated by
\begin{equation}
\int_M\sqrt uu^2\le \Big(\int_Mu^{\frac n4} \Big)^{\frac2n}\Big(\int_Mu^{\frac{2n}{n-2}}\Big)^{\frac{n-2}n},
\end{equation}
where
\begin{equation}
\Big(\int_Mu^{\frac n4} \Big)^{\frac2n}=\Big(\int_M\abs F^\frac{n}2\Big)^\frac2n.
\end{equation}

Applying then the Sobolev inequality
\begin{equation}
\Big(\int_Mu^{\frac{2n}{n-2}}\Big)^{\frac{n-2}n}\le c_1 \int_M\abs{Du}^2+c_2\int_Mu^2,
\end{equation}
\cf \cite{aubin:sobolev},  we obtain
\begin{equation}
\Big(\int_Mu^{\frac{2n}{n-2}}\Big)^{\frac{n-2}n}\le c_3 \Big(\int_M\abs F^\frac{n}2\Big)^\frac2n \Big(\int_Mu^{\frac{2n}{n-2}}\Big)^{\frac{n-2}n},
\end{equation}
where $c_3$ depends on $c_1,c_2,c_0$ and $c$. Hence, we deduce
$u\equiv 0$ or
\begin{equation}
c_3^{-1}\le \Big(\int_M\abs F^\frac{n}2\Big)^\frac2n.
\end{equation}

Setting
\begin{equation}
\ka_0=c_3^{-1}
\end{equation}
finishes the proof.
\ep
\section{The non-compact case}

We now suppose that $M=M^n$ is a complete, non-compact Riemannian manifold. Then there holds
\begin{equation}
H^{1,2}(M)=H^{1,2}_0(M),
\end{equation}
i.e., the test functions $C^\un_c(M)$ are dense in the Sobolev space $H^{1,2}(M)$, see \cite[Lemme 4]{aubin:sobolev} or \cite[Theorem 2.6]{aubin:book}.

Since we do not a priori
\begin{equation}
F^2\in H^{1,2}(M),
\end{equation}
but only
\begin{equation}
F^2\in H^{1,2}_{\tup{loc}}(M),
\end{equation}
the preceding proof has to be modified.

Let $\h=\h(t)$ be defined through
\begin{equation}
\h(t)=
\begin{cases}
1,&t\le 1,\\
(2-t)^q,&1\le t\le 2,\\
0,&t\ge 2,
\end{cases}
\end{equation}
where 
\begin{equation}\lae{3.5}
q=\max(1,\tfrac8n).
\end{equation}

Fix a point $x_0\in M$ and let $r$ be the Riemannian distance function with center in $x_0$
\begin{equation}
r(x)=d(x_0,x).
\end{equation}
Then $r$ is Lipschitz such that
\begin{equation}
\abs{Dr}=1
\end{equation}
almost everywhere.

For $k\ge 1$ define
\begin{equation}
\h_k(x)=\h(k^{-1}r).
\end{equation}
The functions
\begin{equation}
u^{p-1}\h_k^p,
\end{equation}
where
\begin{equation}\lae{3.10}
p=\tfrac n4,
\end{equation}
then have compact support, and multiplying \fre{2.9} with $u^{p-1}\h_k^p$ yields
\begin{equation}
\begin{aligned}
(\tfrac p4+\tfrac18-\e)\int_m\abs{Du}^2u^{p-2}\h_k^p&\le c \Big(\int_M\abs F^\frac n2\Big)^\frac2n \Big(\int_M(u\h_k)^{\frac n{n-2}p}\Big)^\frac {n-2}n\\
&\qq +c_\e\int_M\abs{D\h_k}^2\h_k^{p-2}u^p,
\end{aligned}
\end{equation}
where
$0<\e$ is supposed to be small.

Furthermore, there holds
\begin{equation}
\begin{aligned}
&\int_M\abs{D(u\h_k)^\frac p2}^2=\frac{p^2}4\int_M\abs{Du\h_k+uD\h_k}^2(u\h_k)^{p-2}\\
&\qq\qq\q\le (1+\e)\frac{p^2}4\int_M \abs{Du}^2u^{p-2}\h_k^p+c_\e\frac {p^2}4\int_M\abs{D\h_k}^2\h_k^{p-2}u^p.
\end{aligned}
\end{equation}

Now, choosing $\e$ so small that
\begin{equation}
(1+\e)\tfrac{p^2}4\le p(\tfrac p4+\tfrac18-\e)
\end{equation}
and setting
\begin{equation}
\f=(u\h_k)^\frac p2
\end{equation}
we obtain
\begin{equation}\lae{3.15}
\begin{aligned}
\int_M\abs{D\f}^2\le pc\Big(\int_M\abs F^\frac n2\Big)^\frac2n\Big(\int_M\f^\frac {2n}{n-2}\Big)^\frac {n-2}n+c_\e \int_M\abs{D\h_k}^2\h_k^{p-2}u^p,
\end{aligned}
\end{equation}
where $c_\e$ is a new constant.

We furthermore observe that
\begin{equation}
\abs{D\h_k}^2\h_k^{p-2}\le q^2k^{-2}(2-k^{-1}r)^{qp-2},
\end{equation}
subject to 
\begin{equation}
1\le k^{-1}r\le 2.
\end{equation}

In view of \re{3.5} and \re{3.10}
\begin{equation}
qp-2\ge 0
\end{equation}
and hence
\begin{equation}
\abs{D\h_k}^2\h_k^{p-2}\le q^2k^{-2}.
\end{equation}

Applying now the Sobolev inequality \fre{1.6} to $\f$ and choosing
\begin{equation}
\ka_0=(c_1cp)^{-1}
\end{equation}
we conclude $\abs F\equiv 0$, if
\begin{equation}
\Big(\int_M\abs F^\frac n2\Big)^\frac2n <\ka_0.
\end{equation}

Indeed, if the preceding inequality is valid, then we deduce from \re{3.15}
\begin{equation}
\begin{aligned}
\Big(1-\ka_0^{-1}\Big(\int_M\abs F^\frac n2\Big)^\frac2n\Big)\Big(\int_M\abs\f^\frac{2n}{n-2}\Big)^\frac{n-2}n\le c_\e q^2k^{-2}\int_M\abs F^\frac n2.
\end{aligned}
\end{equation}

In the limit $k\ra\un$ we obtain
\begin{equation}
\Big(\int_M\abs u^\frac{pn}{n-2}\Big)^\frac{n-2}n\le 0.
\end{equation}

\providecommand{\bysame}{\leavevmode\hbox to3em{\hrulefill}\thinspace}
\providecommand{\href}[2]{#2}



\end{document}